\documentclass[11pt,oneside]{article}

\usepackage{graphicx}
\usepackage{amsmath}
\usepackage{amsfonts}
\usepackage{amssymb}
\usepackage{amsthm}
 \usepackage{mathrsfs}
\usepackage{enumerate}
\usepackage{pdfsync}
 \usepackage{stmaryrd}
\usepackage{comment}
\usepackage{color}
\usepackage{bm}
\usepackage{times}
\usepackage{hyperref}
\usepackage{mathtools}
\usepackage[all]{xy}

\usepackage{bm}

\usepackage{tikz}

\bmdefine\taub{\tau}
\bmdefine\mub{\mu}
\bmdefine\lab{\lambda}
\bmdefine\varsigmab{\varsigma}

\numberwithin{equation}{section}
    
\newcommand{\R}{\mathbb{R}}
\newcommand{\E}{\mathbb{E}}

\newcommand{\dd}{\mathrm{d}}
\renewcommand{\P}{\mathbb{P}}

\newcommand{\FF}{\mathcal{F}}
\newcommand{\FFF}{\mathscr{F}}

\newcommand{\diver}{\textnormal{div}}

\newcommand{\vel}{\widetilde{v}}
\newcommand{\velo}{\widehat{v}}

\bmdefine\thetab{\vartheta}
\newcommand{\W}{\mathcal{W}}

\newtheorem{thm}{Theorem}[section]
\newtheorem{lem}[thm]{Lemma}
\newtheorem{prop}[thm]{Proposition}
\newtheorem{cor}[thm]{Corollary}

\theoremstyle{definition}

\theoremstyle{definition}
\newtheorem{assumption}[thm]{Assumption}

\theoremstyle{remark}

\theoremstyle{remark}
\newtheorem{rem}[thm]{Remark}

\vspace{2mm}

\makeatletter
     \newcommand{\xyR}[1]{%
     \makeatletter
     \xydef@\xymatrixrowsep@{#1}
     \makeatother
} 

\makeatletter
     \newcommand{\xyC}[1]{%
     \makeatletter
     \xydef@\xymatrixcolsep@{#1}
     \makeatother
} 

\title{
A trajectorial approach to entropy dissipation \\ for degenerate parabolic equations
}

\vspace{2mm}

\author{  
\textsc{Donghan Kim} \thanks{
Department of Mathematics, University of Michigan (Email: {\it donghank@umich.edu}). 
}  
 \and
\textsc{Lane Chun Yeung}                \thanks{Department of Industrial Engineering \& Operations Research, Columbia University 
(Email: {\it l.yeung@columbia.edu})
}
}

\begin{document}

\maketitle
\bigskip

\begin{abstract} 
\noindent
We consider degenerate diffusion equations of the form $\partial_tp_t = \Delta f(p_t)$ on a bounded domain and subject to no-flux boundary conditions, for a class of nonlinearities $f$ that includes the porous medium equation. We derive for them a trajectorial analogue of the entropy dissipation identity, which describes the rate of entropy dissipation along every path of the diffusion. In line with the recent work \cite{KST20b}, our approach is based on applying stochastic calculus to the underlying probabilistic representations, which in our context are stochastic differential equations with normal reflection on the boundary. This trajectorial approach also leads to a new derivation of the Wasserstein gradient flow property for nonlinear diffusions, as well as to a simple proof of the HWI inequality in the present context. 
\end{abstract}

\smallskip
\noindent{\it Keywords and Phrases:} degenerate diffusion, porous medium equation, entropy dissipation, gradient flow, HWI inequality

\smallskip
\noindent{\it AMS 2000 Subject Classifications: 60H30; 76S05}


\section{Introduction}
\label{sec: intro}

In this paper, we are interested in a class of  quasilinear degenerate  parabolic equations with initial and no-flux boundary conditions of the following form:
\begin{align}	\label{eq: PME2}
	\begin{dcases}
		\partial_t p(t,x) = \Delta \Big(f\big(p(t,x)\big)\Big), \qquad &\text{for } (t, x) \in (0, T) \times U, \\
		~~~ p(0, x) = p_0(x), \qquad \qquad ~~~ &\text{for } x \in \overline{U}, \\
		 \frac{\partial p(t,x)}{\partial n(x)} = 0, \qquad ~~~ &\text{for } (t, x) \in (0, T) \times \partial U,
	\end{dcases}
\end{align}
for a fixed $ T \in (0, \infty) $, an open connected bounded domain $ U \subset \R^d $, and a given initial probability density function $p_0$ on $\overline{U}$. Here, $n(x)$ is the outward normal to the boundary $ \partial U $ at $ x \in \partial U $, and $f : [0, \infty) \rightarrow \mathbb{R}$ is a function representing the nonlinearity. In particular, when $ f(u) = u^m$ for some $ m > 1 $, the partial differential equation of \eqref{eq: PME2} becomes the \emph{porous medium equation}.

Under suitable assumptions on $ f $, it is well known from \cite{carrillo2001entropy} that the solution of \eqref{eq: PME2} converges to a unique stationary distribution, i.e., a probability density function $p_\infty$ satisfying $\Delta \big(f(p_\infty(x))\big) = 0$, and that this convergence can be quantified by the rate of entropy dissipation. More precisely, let us define $ h : (0,\infty) \to \R$ and $ \Phi: [0, \infty) \to \R $ by
\begin{align} \label{eq: h, Phi def}
	h(u) \coloneqq \int_1^u \frac{f'(s)}{s} \, \dd s, \quad \Phi(u) \coloneqq \int_0^u h(s)\,\dd s.
\end{align}
Define also the \emph{entropy functional}
\begin{align}\label{eq: free energy}
\FFF (p) 
& \coloneqq  \int_{U} \Phi \big(p(x) \big)\, \dd x,
\end{align}
for any probability density function $p$ on $\overline{U}$ such that the integral is finite. Then it can be shown that the stationary distribution is the minimizer of  $\mathscr{F}$. Also,  by abbreviating $p_t \coloneqq p(t, \cdot)$, it is well known (see, e.g. \cite[Equation (4)]{carrillo2001entropy}) that
\begin{equation} \label{eq: entro diss iden int} 
	\FFF(p_t) - \FFF(p_{t_0})	
	= -\int_{t_0}^t  I(p_u) \dd u
\end{equation}
holds for every $ 0 \le t_0 \le  t \le T $, where $I$ is the \emph{entropy dissipation functional}, defined by
\begin{align} \label{eq: energy diss func}
    I(p) \coloneqq \int_{U} \left\vert \Phi''\big(p(x)\big)\nabla p(x)\right\vert^2 p(x) \, \dd x,
\end{align}
for any differentiable probability density function $p$ such that the integral is finite. This identity measures the rate of entropy dissipation along the flow of the time-marginal distributions $(p_t)_{t \in [0, T]}$, hence is known as the \emph{entropy dissipation identity}. In particular, the entropy functional $t \mapsto \FFF(p_t)$ is decreasing in time. See also \cite[Lemma 18.14]{Vaz07} and \cite[Equation 3.4]{CT00} for a specific form of this identity for the porous medium equation with drift.

The identity \eqref{eq: entro diss iden int} describes the rate of entropy dissipation at the \emph{ensemble} level of the diffusion modeled by \eqref{eq: PME2}, since it is formulated in terms of the probability distributions $(p_t)$ of the diffusion. The main goal of this paper is to formulate a \emph{trajectorial} analogue to this identity, which describes the rate of entropy dissipation at the level of the individual diffusive particle. To illustrate this, we begin with the following stochastic differential equation (SDE) with normal reflection at the boundary:
\begin{equation}	\label{eq: SDE.intro.0}
	X_t = X_0 + \int_{0}^t\sqrt{\frac{2f\big(p(s, X_s)\big)}{p(s, X_s)}}\dd W_s -  \int_{0}^t n(X_s)\, \dd L_s \in \overline{U},   \qquad X_0 \sim p_0.
\end{equation}
Here, $W$ is a $d$-dimensional standard Brownian motion and $ L $ is a nondecreasing continuous process satisfying
\begin{align}
	L_t = \int_0^t \mathrm{1}_{\{X_s \in \partial \Omega\}} \, \dd L_s, \quad L_0 = 0.
\end{align}
The stochastic process $(X_t)_{0 \le t \le T}$ is probabilistic representation of \eqref{eq: PME2}, in the sense that its time-marginal probability density function is given by the solution $p(t, \cdot)$ of \eqref{eq: PME2}. Intuitively, the diffusion $X_t$ is reflected on the boundary $\partial U$ in the direction $-n(X_t)$. The reflecting term $L$ is associated with a multidimensional analogue of the \emph{local time} on $\partial U$ \cite{sato1962local}. With this probabilistic representation, the entropy at time $t$ can then be expressed as an expectation
\begin{equation}\label{eq: free energy 3}
\FFF (p_t) = \int_{U} v(t, x)  \, p(t, x)  \, \mathrm{d}x =\E \big[ v(t, X_t) \big], \quad \text{where}\footnote{In the case of the porous medium equation, $ v $ is known as the \emph{pressure function}.}  v(t,x) \coloneqq \frac{\Phi\big(p(t,x)\big)}{p(t,x)}.
\end{equation}
Using stochastic calculus, we shall derive the dynamics of the \emph{entropy  process} $\big(v(t, X_t)\big)_{t\in [0,T]}$, in terms of the \emph{semimartingale decomposition}
	\begin{equation}  \label{eq: intro semimart decom}
		v(t, X_t)  - v(0, X_{0}) = M_{t} + F_{t}, \qquad \textnormal{for} \quad 0 \le t \le T, 
	\end{equation}
where $M$ is a martingale and $F$ is a process of finite variation. This decomposition describes the evolution of the entropy process along every trajectory of the particle, thus it can be seen as a trajectorial analogue of \eqref{eq: entro diss iden int}. In fact, \eqref{eq: entro diss iden int} can be recovered from \eqref{eq: intro semimart decom} by averaging over these trajectories; in other words, by taking expectation.

Our work is much inspired from the recent work \cite{KST20b}, which provides a trajectorial approach to the relative entropy dissipation for Fokker-Planck equations.
Subsequently, this approach has been extended to Markov chains \cite{KMS20} and to McKean-Vlasov equations \cite{tschiderer2021trajectorial}. It is therefore natural to expect an adaptation of the approach for the porous media type equation \eqref{eq: PME2}. Compared with prior work, a key difficulty in our setting stems from the degenerate parabolicity of \eqref{eq: PME2}. More specifically, as $ f' $ is not assumed to be bounded from below by some strictly positive constant, \eqref{eq: PME2} is not uniformly parabolic. Without uniform parabolicity, equations of this form are only expected to have weak solutions \cite{wilhelm1983quasilinear}, but not classical solutions. However, such regularity is crucial for applying It\^{o} calculus. To this end, we require the initial condition $ p_0 $ to be nondegenerate, i.e., $ \kappa^{-1} \le p_0(x)\le \kappa  $ for some $ \kappa > 1 $. This will ensure that \eqref{eq: PME2} has a smooth solution. Also, besides considering diffusions on a bounded domain, another main difference with the prior work \cite{KST20b,tschiderer2021trajectorial} is that our main trajectorial result (Theorem \ref{thm: traj} below) is stated in the forward direction of time.

Along with our trajectorial approach come two applications. The first application is a new derivation of the Wasserstein gradient flow property of \eqref{eq: PME2}, which states that the curve of time-marginal probability density functions $ (p_t)_{t \in [0,T]} $ of \eqref{eq: PME2} descends in the steepest possible direction of the entropy functional $ \mathscr{F} $ in $\mathcal{P}(\overline{U})$, the space of probability measures on $U$. Here, $\mathcal{P}(\overline{U})$ is equipped with the quadratic Wasserstein distance $ \mathcal{W}_2 $, defined by
\begin{equation} \label{eq: wasserstein dist}
	\W_{2}(\mu,\nu) \coloneqq \sqrt{ \inf_{\pi } \int_{U \times U } \vert x-y \vert^{2} \, \pi(\dd x,\dd y) }, \quad \text{for any }\mu, \nu \in \mathcal{P}(\overline{U}),
\end{equation}
where the infimum is taken over $ \pi \in \mathcal{P}(\overline{U} \times \overline{U})$ with marginals $ \mu $ and $ \nu $.

For the porous medium equation on  $\mathbb{R}^d$, this property was discovered by Otto in his seminal paper \cite{Ott01}, where he introduced a formal Riemannian structure on $\mathcal{P}(\R^d)$. More recently, Ambrosio, Gigli and Savar\'e \cite{AGS08} developed a rigorous theory of gradient flows on general metric spaces based on the notion of curves of maximal slopes. Similar results have been established for porous medium equations on discrete spaces \cite{EM14} and with fractional pressure \cite{lisini2018gradient}.  

To show the gradient flow property, we adopt the methodology in \cite{KST20b} of \emph{perturbing} the SDE \eqref{eq: SDE.intro.0} from some time $t_0 \in [0, T)$ by adding a gradient drift $ \nabla \beta $:
 \begin{align*}
 	X^{\beta}_t= X^{\beta}_{t_0 } - \int_{t_0}^t \nabla \beta(X^{\beta}_s) \, \dd s + \int_{t_0}^t \sqrt{\frac{2f\big(p^{\beta}(s, X^{\beta}_s)\big)}{p^{\beta}(s, X^{\beta}_s)}} \, \dd W^{\beta}_{s} -  \int_{t_0}^t n(X^{\beta}_s)\, \dd L^\beta_s.
 \end{align*}
Here, $p^{\beta}(t, \cdot)$ is the time-marginal probability density function for the solution $(X^{\beta}_t)_{t_0 \le t \le T}$ of this perturbed SDE. By deriving the dynamics of the associated perturbed entropy process, we obtain an analogous entropy dissipation identity for the perturbed diffusion. On the other hand, we can also explicitly compute the rates of changes of the Wasserstein distances along both the perturbed curve $(p^\beta_t)$ and the unperturbed curve $(p_t)$. Thus, the entropy dissipation rates can be measured not in terms of time elapsed, but in terms of the Wasserstein distances traveled by the curve of time-marginal probability density functions, both in the perturbed and unperturbed settings. Comparing them allows us to establish the maximal rate of entropy dissipation for the unperturbed diffusion \eqref{eq: PME2}, by measuring the exact effect of each perturbation.

The second application of the trajectorial approach is a simple proof of the HWI inequality in the context of the nonlinear equation \eqref{eq: PME2}, which is a special case of \cite[Theorem 4.2]{AGK04}.  It is an interpolation inequality relating the entropy functional (H)\footnote{The letter ``H" comes from the choice $f(u) = u$ in \eqref{eq: h, Phi def}, in which case the entropy functional defined in \eqref{eq: free energy}  satisfies $\mathscr{F}(p) = H(p) - 1$, where $H(p) = \int_U p(x) \log p(x) \, \dd x$ is the (negative of the) differential entropy.}, the Wasserstein distance (W) and the entropy dissipation functional (I). More precisely, it states that  
\begin{align} \label{eq: HWI}
	\FFF(\rho_0) - \FFF(\rho_1) 
	& \le \sqrt{I(\rho_0)} \, \W_2(\rho_0, \rho_1)
\end{align}
holds for any $ \rho_0, \rho_1 \in \mathcal{P}(\overline{U})$.
We will prove this inequality by applying a trajectorial approach similar to the one just described, but instead to the \emph{displacement interpolation} between $ \rho_0$ and $ \rho_1 $.

The rest of the paper is organized as follows. In Section \ref{subsec: stochastic representation}, we introduce our setup and state some preliminary lemmas. Our main result is stated in Section \ref{subsec: traj PME}.  Section \ref{subsec: gradient flow} formulates the gradient flow property, while Section \ref{subsec: HWI} develops the HWI inequality. Proofs are provided in Section \ref{sec: proofs}. 

\medskip

\section{Setting and main results} \label{sec: setup}

\subsection{Setup}	\label{subsec: stochastic representation}

We impose the following assumption on the initial distribution $p_0$ and the nonlinearity $ f $ of the degenerate parabolic equation~\eqref{eq: PME2}. 

\begin{assumption}	\label{assum1}
	\mbox{}
	\begin{enumerate}
		\item [(a)]  The domain $ U $ is an open connected bounded subset of $ \R^d $ for some $ d \ge 2 $, and the boundary  $ \partial U$ is smooth.
		\item [(b)]	The initial datum $p_0$ is a smooth probability density function on $\overline{U}$. Moreover, it is non-degenerate, i.e., there exists $ \kappa > 1$ for which $ \kappa^{-1} \le p_0(x)\le \kappa  $ holds for all $ x \in \overline{U}$. Also,
		\begin{align} \label{eq: assum noflux initial}
			\frac{\partial p_0(x)}{\partial n(x)} = 0, \quad \text{for all } x \in \partial U.
		\end{align}
	    where $ n(x) $ is the outward normal to the boundary $ \partial U $.

		\item [(c)] The function $ f : [0, \infty) \to \R $ is smooth and strictly increasing. Moreover, $ f(0) = f'(0) = 0 $ and $ f'(u) > 0 $ for all $ u > 0 $. Also, $ f' $ is nondecreasing. 
		\item[(d)] The function $h$, defined in \eqref{eq: h, Phi def}, belongs to $L^1_{\text{loc}}([0, \infty))$.
		
	\end{enumerate}
\end{assumption}

\begin{rem}
	Assumption \ref{assum1}(c) -- (d) covers the porous medium equations, in which case $ f(u) = u^{m} $ for some $ m > 1 $. The assumption $f'(u) > 0$ implies that \eqref{eq: PME2} is a \emph{parabolic} partial differential equation (PDE), while the assumption $ f'(0) = 0 $ implies that \eqref{eq: PME2} is \emph{degenerate} parabolic. 
\end{rem}
We first collect some basic properties of the solution to the PDE~\eqref{eq: PME2} in the following lemma. These properties are classical, and we refer to Chapter 3 of the  monograph  \cite{Vaz07} and the references therein for a comprehensive overview.

\begin{lem} \label{lem: sol of PME}
	Under Assumption~\ref{assum1}, there exists a smooth solution $p \in C^{\infty}([0, T] \times \overline{U})$ of \eqref{eq: PME2}. Moreover, $\int_U p(t,x)\, \dd x = 1$ for all $t \in [0,T]$ and $\kappa^{-1} \le p(t, x)\le \kappa$ for all $ (t,x) \in [0, T] \times \overline{U} $.
\end{lem}

For the rest of the paper, we fix $p$ as the solution given in Lemma \ref{lem: sol of PME}. Fix also a filtered probability space $\big(\Omega,\mathbb{F},  \mathcal{F} = (\mathcal{F}_t)_{0 \le t \le T}, \P \big)$ supporting a $\mathcal{F}$-Brownian motion $W$ and a $\mathcal{F}_0$-measurable random vector  $\xi : \Omega \to \R^d$ with $\P \circ \xi^{-1} = p_0$. We shall denote by $\E$ the expectation taken with respect to $\P$. We will make use of the probabilistic representation of \eqref{eq: PME2}, which is described in Lemma \ref{lem: stoch.rep} below in terms of the solution of the following SDE with normal reflection on the boundary:
\begin{equation}	\label{eq: SDE.intro}
	\begin{dcases}
		X_t = X_0 + \int_{0}^t\sqrt{\frac{2f\big(p(s, X_s)\big)}{p(s, X_s)}} \, \dd W_s -  \int_{0}^t n(X_s)\, \dd L_s \in \overline{U},  \qquad t \in [0, T],\\
		X_0 = \xi, \\
		L_t = \int_0^t \mathrm{1}_{\{X_s \in \partial U\}} \, \dd L_s, \quad
	    [0, T] \ni t \to L_t \text{ is nondecreasing,  continuous with } L_0 = 0.
	\end{dcases}
\end{equation}
For an introduction to SDEs with reflection and their connections with nonlinear paraoblic PDEs, we refer to the lecture notes \cite{pilipenko2014introduction} and \cite{soner2007stochastic}.

The following lemma shows that \eqref{eq: SDE.intro} is well-posed and provides the probabilistic representation of  \eqref{eq: PME2}. Similar results are known for the porous medium equations \cite{benachour1996process} as well as for general nonlinear equations of the form \eqref{eq: PME2} with discontinuous coefficients  \cite{barbu2011probabilistic, blanchard2010probabilistic}, or with the half-line as the domain \cite{ciotir2014probabilistic}. See also \cite{hu2017gradient} for a martingale method for establishing gradient estimates for the porous medium and fast diffusion equations.

\begin{lem} \label{lem: stoch.rep}
	Suppose Assumption \ref{assum1} holds. Then the SDE with reflection \eqref{eq: SDE.intro} has a pathwise unique, strong solution $ (X,L) $, for which the probability density functions of $ (X_t)_{t \in [0,T]} $  are given by $\big(p(t, \cdot)\big)_{t \in [0,T]}$.
\end{lem}

\medskip

\subsection{Trajectorial entropy dissipation of degenerate parabolic equation}  \label{subsec: traj PME}

Our first main result is the dynamics of entropy dissipation along every trajectory of the diffusion \eqref{eq: SDE.intro}, formulated in terms of the semimartingale decomposition of the entropy process $ \big(v(t, X_t)\big)_{t \in [0, T]}$.

To begin, we define the \emph{entropy dissipation function}
\begin{equation} \label{eq: D process}
	D(t, x) := \bigg(\varphi'(p) \Delta f(p) +  \frac{f(p)}{p} \Delta v\bigg) (t, x), \quad (t,x)\in [0, T] \times \overline{U},  \quad \text{where} \quad \varphi (u) := \frac{\Phi(u)}{u},
\end{equation}
and $\Phi$ is defined in \eqref{eq: h, Phi def}. This function will help provide the exact trajectorial rate of entropy dissipation, as will be explained in Remark \ref{rem: cond traj} below.

\begin{thm}  \label{thm: traj}
	Suppose Assumption~\ref{assum1} holds. Then the entropy process $ \big(v(t, X_t)\big)_{t \in [0, T]}$ admits the semimartingale decomposition 
	\begin{equation} \label{eq: semi-mart decom}
		v(t, X_t)  - v(0, X_0) = M_t + F_t, \qquad \textnormal{for} \quad t \in [0,T]. 
	\end{equation}
where 
\begin{align}	\label{eq: cum en diss process}
	F_t \coloneqq  \int_0^t D (s,X_s)\,  \dd s, \quad \text{and} \quad M_t \coloneqq  \, \int_0^t \Bigg\langle \sqrt{\frac{2f\big(p(s, X_s)\big)}{p(s, X_s)}}\nabla v(s,X_s), \, \dd W_s \Bigg\rangle
\end{align}
is an $L^2$- bounded martingale. Also, we have
\begin{equation} \label{eq: fish.inf.exp}
	\E \big[F_{t}\big] 
	= -\int_{0}^{t} I\big( p_t\big)  \, \dd t 
	> -\infty \, , \qquad \textnormal{for} \quad t \in [0,T].
\end{equation}
\end{thm}
The proof of Theorem~\ref{thm: traj} will be given in Section \ref{subsec: proof.traj perturbed}. This result is the analogue of  \cite[Theorem 4.1]{KST20b} and \cite[Theorem 3.1]{tschiderer2021trajectorial}. In contrast to them, Theorem \ref{thm: traj} here is stated in the forward direction of time. 

By aggregating this trajectorial result, i.e., taking expectation, we recover the entropy dissipation identity \eqref{eq: entro diss iden int} and its differential version \eqref{eq: entro diss iden2}. Furthermore, by taking \emph{conditional} expectation, we obtain below the conditional trajectorial rate of entropy dissipation \eqref{eq: conditional entro diss}.

\begin{cor} \label{cor: deBruijn}
	Suppose Assumption~\ref{assum1} holds. Then for every $0 \le t_0 \le t \le T$,  the entropy dissipation identity \eqref{eq: entro diss iden int} holds. The corresponding differential version
	\begin{equation} \label{eq: entro diss iden2}
		\frac{\mathrm{d}}{\mathrm{d}t} \bigg\vert_{t=t_0}\FFF(p_t) = - I(p_{t_0})
	\end{equation}
also holds for all $t_0 \in [0, T]$. Moreover, for all $t_0 \in [0, T]$, the conditional trajectorial rate of entropy dissipation is given by
\begin{equation}    \label{eq: conditional entro diss}
	\lim_{t \downarrow t_0} \frac{\mathbb{E} [ v(t, X_{t}) \, \big| \, \FF_{t_0} ] - v(t_0, X_{t_0})}{t-t_0}
	= D(t_0, X_{t_0}),
\end{equation}
where the limit exists in $ L^1(\P) $.
\end{cor}
\begin{rem} \label{rem: cond traj}
From \eqref{eq: fish.inf.exp}, we see that $\E[D(t_0, X_{t_0})] = -I(p_{t_0})$ holds. This explains why \eqref{eq: conditional entro diss} constitutes a conditional trajectorial version of the entropy dissipation identity \eqref{eq: entro diss iden2}. 
\end{rem}

\medskip

\subsection{Gradient flow property of the degenerate parabolic equation, via perturbation analysis} \label{subsec: gradient flow}
In this subsection, we discuss how our trajectorial approach leads to a new interpretation of the Wasserstein gradient flow property of the degenerate parabolic equation \eqref{eq: PME2}. Following the method of \cite{KST20b,tschiderer2021trajectorial}, we shall \emph{perturb} the degenerate parabolic equation. 

To this effect, let $\beta$ be a perturbation potential satisfying the following assumption.

\begin{assumption}	\label{assum.perturbation}
	The perturbation potential $\beta : \overline{U} \to \R$ is smooth. Moreover, $ \nabla \beta(x) = 0 $ for $ x \in \partial U $.
\end{assumption}

For the rest of the paper, we fix a $t_0 \in [0, T)$ and a perturbation $\beta$ satisfying Assumption~\ref{assum.perturbation}. Consider the following Neumann problem, which can be viewed as a perturbed version of \eqref{eq: PME2}:
\begin{align}	\label{eq: PME perturbed}
	\begin{dcases}
		\partial_t p^{\beta}(t,x) 
		= \diver \Big( \nabla f\big(p^{\beta}(t,x)\big) + p^{\beta}(t,x)  \nabla\beta(x)  \Big), \quad &\text{for } (t, x) \in (t_0, T] \times U, \\
		p^{\beta}(t_0, x) = p(t_0, x),  \quad &\text{for } x \in \overline{U},\\
		\frac{\partial p^{\beta}(t,x)}{\partial n(x)} = 0, \quad &\text{for } x \in \partial U.
	\end{dcases}
\end{align}
The following result, which is the ``perturbed analogue'' of Lemma \ref{lem: sol of PME}, shows the existence of a strictly positive smooth solution to \eqref{eq: PME perturbed} in a short time interval.

\begin{lem} \label{lem: sol of p. PME}
	Under Assumptions \ref{assum1} and \ref{assum.perturbation}, there exists $T_\beta \in (t_0, T]$ such that  \eqref{eq: PME perturbed} has a smooth solution in $p^\beta \in C^{\infty}\big( [t_0, T_\beta] \times \overline{U} \big)$. Moreover, $\int_U p^\beta(t,x)\, \dd x = 1$ for all $t \in [t_0, T_\beta]$ and $ \frac{1}{2\kappa}\le p^\beta(t, x)\le \frac{3}{2\kappa}$ for all $ (t, x) \in [t_0, T_\beta] \times \overline{U} $. 
\end{lem}

For the rest of the paper, we fix $ p^\beta $ as  given by Lemma \ref{lem: sol of p. PME}. The corresponding probabilistic representation of the perturbed PDE \eqref{eq: PME perturbed} is the following SDE with reflection:
\begin{equation}	\label{eq: SDE.p}
	\begin{dcases}
		X^\beta_t= X^\beta_{t_0 } - \int_{t_0}^t \nabla \beta(X^\beta_s) \, \dd s + \int_{t_0}^t \sqrt{\frac{2f\big(p^{\beta}(s, X^\beta_s)\big)}{p^{\beta}(s, X^\beta_s)}} \, \dd W_{s} -  \int_{t_0}^t n(X^\beta_s)\, \dd L^\beta_s \in \overline{U} , \quad t \in [t_0, T_{\beta}],\\ 
		X^\beta_{t_0} = X_{t_0},  \\
		L^\beta_t = \int_{t_0}^t \mathrm{1}_{\{X^\beta_s \in \partial U\}} \, \dd L^\beta_s, \quad
		[t_0, T_{\beta}] \ni t \mapsto  L^\beta_{t} \text{ is nondecreasing continuous with } L^\beta_{t_0} = 0.
	\end{dcases}
\end{equation}

By analogy with Lemma \ref{lem: stoch.rep}, the following result ensures that \eqref{eq: SDE.p} is well-posed and is the stochastic counterpart of \eqref{eq: PME perturbed}.
\begin{lem} \label{lem: stoch.rep.p}
	Under Assumptions \ref{assum1} and \ref{assum.perturbation}, the SDE with reflection \eqref{eq: SDE.p} has a pathwise unique, strong solution $ (X^\beta_t,L^\beta_t)_{t \in [t_0,T_\beta]} $, for which the probability density functions of $ (X^\beta_t)_{t \in [t_0, T_\beta]}$  are given by $\big(p^\beta(t, \cdot)\big)_{t \in [t_0,T_\beta]}$.
\end{lem}

As before, let us abbreviate $p^\beta_t \coloneqq p^\beta(t, \cdot)$.  In parallel to \eqref{eq: free energy 3}, we can express the entropy for the perturbed diffusion at time $t$ as
\begin{align} \label{def: pressure func.p}
\FFF \big(p^{\beta}_t\big) =  \int_{U} v^{\beta}(t, x) \, p^{\beta}(t, x) \, \mathrm{d}x =  \E \big[ v^{\beta}(t, X^\beta_t) \big], \quad \text{where} \quad 	v^\beta(t, x) \coloneqq \frac{\Phi \big(p^{\beta}(t, x) \big)}{p^{\beta}(t, x)}.
\end{align}
Therefore, we similarly call $ \big(v^\beta(t, X^\beta_t)\big)_{t \in [t_0, T_\beta]}$ the \emph{perturbed entropy process}. The following result, which is the perturbed counterpart of Theorem \ref{thm: traj}, derives the dynamics of this process. By analogy with \eqref{eq: D process}, we introduce the \emph{perturbed entropy dissipation function} 
\begin{align} \label{eq: en diss process.p.r}
	D^{\beta}(t, x) \coloneqq & \bigg(\varphi'\big(p^{\beta}\big) \, \diver \Big( \nabla f\big( p^{\beta}\big) + p^{\beta} \nabla \beta(x) \Big)  + \frac{f\big(p^{\beta}\big)}{p^{\beta}} \Delta v^{\beta} - \langle \nabla v^{\beta}, \, \nabla \beta\rangle\bigg)(t, x),
\end{align}
for $ (t,x) \in [t_0, T_\beta] \times \overline{U}$. 

\begin{thm}  \label{thm: traj.p.r.}
	Suppose Assumptions~\ref{assum1} and \ref{assum.perturbation} hold. Then the perturbed entropy process $\big(v^{\beta}(t, X^{\beta}_t)\big)_{t \in [t_0, T_\beta]}$ admits the semimartingale decomposition 
	\begin{equation} \label{eq: semi-mart decom.p}
		v^{\beta}(t, X^{\beta}_t)  - v^{\beta}(t_0, X^{\beta}_{t_0}) = M^{\beta}_t + F^{\beta}_t, \qquad \text{for} \quad t \in [t_0, T_{\beta}],
	\end{equation}
	where 
	\begin{equation} \label{eq: cum en diss process.p}
	F^{\beta}_{t} \coloneqq \int_{t_0}^t {D}^{\beta}(s, \, X^{\beta}_s) \,  \dd s, \quad {and} \quad 	M^{\beta}_t \coloneqq \int_{t_0}^t \Bigg\langle \sqrt{\frac{2f\big(p^{\beta}(s, X^{\beta}_s)\big)}{p^{\beta}(s, X^{\beta}_s)}} \nabla v^{\beta} (s, X^{\beta}_s), \, \dd W_s  \Bigg\rangle
	\end{equation}
	is an $L^2$-bounded martingale. Also, we have for $t \in [t_0,  T_{\beta}]$,
	\begin{equation} \label{eq: fish.inf.exp.p}
		\E \big[F^{\beta}_{t}\big] 
		=  -\int_{t_0}^{t} I( p^{\beta}_s)  \, \dd s  - \int_{t_0}^t  \E \Big[ \big\langle \nabla h \big(p^{\beta}(s, X^{\beta}_s) \big), \, \nabla \beta(X^{\beta}_s) \big\rangle \Big] \, \dd s
		> - \infty.
	\end{equation}
\end{thm}

Once again, by averaging this trajectorial result, we obtain the following perturbed entropy dissipation identity and its conditional trajectorial version, to which a comment similar to Remark \ref{rem: cond traj} applies.

\begin{cor} \label{cor: deBruijn.p}
	Suppose Assumptions~\ref{assum1} and \ref{assum.perturbation} hold. For every $ t \in [t_0, T_{\beta}]$, the following perturbed entropy dissipation identity holds:
	\begin{equation}			\label{eq: entro diss integral.p}
		\FFF\big(p^{\beta}_t\big) - \FFF\big(p^{\beta}_{t_0}\big)	
		= -\int_{t_0}^t  I\big(p^{\beta}_s\big) \mathrm{d} s
		-  \int_{t_0}^t  \E \Big[ \big\langle \nabla h\big(p^{\beta}(s, X^{\beta}_s)\big), \, \nabla \beta(X^{\beta}_s) \big\rangle \Big] \, \dd s. 
	\end{equation}
The corresponding differential version also holds:
	\begin{align} 
		\frac{\mathrm{d}}{\mathrm{d}t} \bigg\vert_{t=t_0^+} \FFF(p^{\beta}_t) 
		& = - I(p_{t_0}) - \E \Big[ \big\langle \nabla h \big(p(t_0, X_{t_0}) \big), \, \nabla \beta(X_{t_0}) \big\rangle \Big].      \label{eq: entro diss iden.p}
	\end{align}
Moreover, the conditional trajectorial rate of entropy dissipation for the perturbed diffusion is given by
\begin{equation}    \label{eq: conditional entro diss perturb 1}
	\lim_{t \downarrow t_0} \frac{\mathbb{E} \Big[ v^{\beta}(t, X^{\beta}_{t}) \, \big| \, \FF_{t_0} \Big] - v^{\beta}(t_0, X_{t_0})}{t-t_0}
	= D^{\beta}(t_0, X_{t_0}),
\end{equation}
where the limit exists in $L^1$.
\end{cor}

The last ingredients we need are the rates of change of the Wasserstein distances along the curve of the marginal distributions $ (p_t) $ and $ (p^\beta_t) $. The following result is a consequence of the general theory of Wasserstein metric derivatives for absolutely  continuous curves \cite[Chapter 8]{AGS08}, but our setting allows for a more direct proof, which we provide in Section \ref{sec: Wasserstein}. We recall the definitions of $h$ in \eqref{eq: h, Phi def} and the perturbation potential $\beta$ described at the beginning of Section \ref{subsec: gradient flow}.
\begin{lem} \label{thm: wass}
	Suppose Assumptions~\ref{assum1} and \ref{assum.perturbation} hold. Then the Wasserstein metric slope along the unperturbed curve $(p_t)$ is given by
\begin{equation} \label{eq: wass.limit}
	\lim_{t \downarrow t_0} \frac{\W_2(p_t, p_{t_0})}{t - t_0} =  \Big\Vert \nabla h\big(p(t_0, X_{t_0})\big) \Big\Vert_{L^2}.
\end{equation}
Similarly, the Wasserstein metric slope along the perturbed curve $\big(p^\beta_t\big)$ is given by
\begin{equation} \label{eq: wass.limit.perturb}
	\lim_{t \downarrow t_0} \frac{\W_2 \big(p^{\beta}_t, p^{\beta}_{t_0}\big)}{t - t_0} =  \Big\Vert \nabla h \big(p(t_0, X_{t_0})\big)  + \nabla  \beta(X_{t_0}) \Big\Vert_{L^2}.
\end{equation}
\end{lem}

Combining the entropy dissipation identities \eqref{eq: entro diss iden2} and \eqref{eq: entro diss iden.p} with the Wasserstein derivatives \eqref{eq: wass.limit}--\eqref{eq: wass.limit.perturb} allows us to derive the gradient flow property. We define the Wasserstein metric slopes of the entropy functional $ \mathscr{F} $ along the unperturbed curve $ (p_t) $ and along the perturbed curve $ (p^\beta_t) $, respectively, by
\begin{align}
	|\partial \mathscr{F}|_{\W_2} (p_{t_0}) \coloneqq \lim_{t \downarrow t_0}   \frac{\FFF(p_t) - \FFF(p_{t_0})}{\W_2 (p_t, p_{t_0})} \quad \text{and} \quad  |\partial \mathscr{F}|_{\W_2} (p^\beta_{t_0}) \coloneqq   \lim_{t \downarrow t_0} \frac{\FFF\big(p^{\beta}_t\big) - \FFF \big(p^{\beta}_{t_0}\big)}{\W_2 (p^{\beta}_t, p^{\beta}_{t_0})}. 
\end{align}
The following theorem computes both of these slopes explicitly, which shows in particular that the unperturbed slope $ |\partial \mathscr{F}|_{\W_2} (p_{t_0})  $ is always steeper than the perturbed slope $  |\partial \mathscr{F}|_{\W_2} (p^\beta_{t_0}) $.
\begin{thm}\label{thm: steepest.decsent}
    Suppose Assumptions~\ref{assum1} and \ref{assum.perturbation} hold. Then the Wasserstein metric slope of the entropy functional along the unperturbed curve $(p_t)$ is given by
    \begin{align}
        |\partial \mathscr{F}|_{\W_2} (p_{t_0}) 
	    &=  -\Big\Vert \nabla h\big(p(t_0, X_{t_0})\big) \Big\Vert_{L^2}= -\sqrt{I(p_{t_0})}.									\label{eq: FW}
    \end{align}
Similarly, if $||\nabla h\big(p(t_0, X_{t_0})\big) + \nabla \beta(X_{t_0})||_{L^2} > 0$, then the Wasserstein metric slope along the perturbed curve $\big(p^\beta_t\big)$ is given by
\begin{align}
        |\partial \mathscr{F}|_{\W_2} \big(p^\beta_{t_0}\big)
	    &=  - \bigg\langle \nabla h\big(p(t_0, X_{t_0})\big), \, \frac{\nabla h\big(p(t_0, X_{t_0})\big) + \nabla \beta(X_{t_0})}{\Big\Vert\nabla h\big(p(t_0, X_{t_0})\big) + \nabla \beta(X_{t_0}) \Big\Vert_{L^2}} \bigg\rangle_{L^2}.	\label{eq: FW.p}
    \end{align}
In particular, 
    \begin{align}
	    	|\partial \mathscr{F}|_{\W_2} (p_{t_0}) \le |\partial \mathscr{F}|_{\W_2} \big(p^\beta_{t_0}\big), \label{eq: FW-FW.p} 
    \end{align} 
and equality holds if and only if $\nabla h\big(p(t_0, X_{t_0})\big) + \nabla \beta(X_{t_0})$ is a.s. a scalar multiple of $\nabla h\big(p(t_0, X_{t_0})\big)$.
\end{thm}

\medskip

\subsection{HWI inequality} \label{subsec: HWI}
In this subsection, we apply a similar trajectorial approach to derive the HWI inequality \eqref{eq: HWI}. Let us fix two probability density functions $\rho_0, \rho_1 \in \mathcal{P}(\overline{U})$ and impose the following assumption.
\begin{assumption} \label{FWI assum} 
    \mbox{}
	\begin{enumerate}
		\item [(a)] Both $\rho_0, \rho_1$ are strictly positive and smooth. Also, $ \rho_0(x) = \rho_1(x)$ for all $ x \in \partial U $.
		\item [(b)] The function $ f : [0, \infty) \to \R $ is smooth. Moreover, the function $h$, defined in \eqref{eq: h, Phi def}, belongs to $L^1_{\text{loc}}([0, \infty))$.
		\item [(c)] The function $ r \mapsto r^d \Phi (r^{-d})$ is convex nonincreasing on $ (0, \infty) $, where $\Phi$ is defined in  \eqref{eq: h, Phi def}. 
	\end{enumerate}
\end{assumption}
\begin{rem}
Assumption (b), (c) is satisfied by the porous medium equation, see \cite[Examples 5.19]{Vil03}.
\end{rem}
By Brenier's theorem \cite[Theorem 2.12(ii)]{Vil03}, there exists a convex function $\psi: \overline{U} \to \R$ such that $ \nabla \psi $ is the optimal transport map from $ \rho_0 $ to $ \rho_1$, i.e.,
\begin{equation}    \label{eq: W with varphi}
	\W_2^2(\rho_0, \rho_1) = \int_{U} \big\vert x - \nabla \psi (x) \big\vert^2 \, \rho_0(\dd x).
\end{equation}
Let $ (\rho_t)_{t \in (0,1)} $  denote the displacement interpolation between $ \rho_0 $ and $ \rho_1 $, i.e.,
\begin{align*}
	\rho_t = \rho_0  \circ \big( (1 - t) \text{Id} + t \nabla \psi \big)^{-1}, \quad \text{for } t \in (0,1).
\end{align*}
It is known that each $\rho_t$ has a probability density function \cite[Remarks 5.13(i)]{Vil03}. For the following result, we recall the entropy functional $\mathscr{F}$ defined in \eqref{eq: free energy}, the nonlinearity $f$ satisfying Assumption \ref{FWI assum}(b), and the convex function $\psi$ described just above.

\begin{prop} \label{prop: FWI} Suppose Assumption \ref{assum1}(a) and Assumption \ref{FWI assum}(a--b) hold. Then the rate of change of $t \mapsto \FFF(\rho_t)$ at $t=0$ is given by
	\begin{equation} \label{eq: rate of change F}
		\frac{\dd}{\dd t} \Big|^{+}_{t =0} \FFF(\rho_t) 
		= \int_{U} \big\langle \nabla f\big(\rho_0(z)\big), \nabla \psi(z) - z \big\rangle \, \dd z.
	\end{equation}
\end{prop}
Using this proposition and the displacement convexity of the entropy functional $ \mathscr{F} $, we obtain the HWI inequality \eqref{eq: HWI}. 
\begin{thm} \label{thm: FWI}
	Suppose Assumption \ref{assum1}(a) and Assumption \ref{FWI assum} hold. Then
	\begin{align} \label{eq: first.func.ineq}
		\FFF(\rho_0) - \FFF(\rho_1) &\le -\int_{U} \big\langle \nabla f\big(\rho_0(z)\big), \nabla \psi(z) - z \big\rangle \, \dd z \le \sqrt{I(\rho_0)} \, \W_2(\rho_0, \rho_1),
	\end{align}
	where $I$ is the entropy dissipation functional defined in \eqref{eq: energy diss func}.
\end{thm}

\medskip

\section{Proofs} \label{sec: proofs}

\subsection{Proofs of Lemmas \ref{lem: sol of PME}  and \ref{lem: sol of p. PME}}

\smallskip

{\noindent \textbf{Proof of Lemma \ref{lem: sol of PME}}.}
We adopt the same method as in the proof of \cite[Theorem 3.1]{Vaz07}, which exploits the nondegeneracy of the initial condition of Assumption \ref{assum1}(b) in a crucial manner . Let  $ \tilde{f} : [0, \infty) \to \R $ be a smooth function satisfying $ \tilde{f}(u) = f(u)$ for $ \kappa^{-1} \le u \le \kappa$, $ \tilde{f}'(u) > \epsilon^{-1} $ for $ 0 \le u \le \kappa^{-1} $ and $ \tilde{f}'(u) < \epsilon $ for $ u \ge \kappa $, where $ \epsilon > 1 $ is some fixed constant. Consider the PDE 
	\begin{align} \label{eq: PDE aux}
		\partial_t p(t,x) = \Delta \Big(\tilde{f}\big(p(t,x)\big)\Big), \qquad \text{for } (t, x) \in (0, T] \times U
	\end{align}
subject to the same initial and Neumann boundary conditions as in \eqref{eq: PME2}. By Assumption \ref{assum1}(c), $ f' $ is increasing, so $ \tilde{f}'(u) > \epsilon^{-1}$ for all $ u \ge 0 $. This implies that \eqref{eq: PDE aux} is uniformly parabolic. We can therefore apply standard quasilinear theory \cite[Chapter 3.1]{Vaz07} to obtain a smooth solution $p \in C^{\infty}([0,T] \times \overline{U})$ to \eqref{eq: PDE aux}. By the comparison principle, $\kappa^{-1} \le p \le \kappa$, so $ p $ also satisfies \eqref{eq: PME2}. 

Finally, the mass conservation law \cite[Chapter 3.3.3]{Vaz07} implies that the total mass $\int_U p(t,x)\, \dd x = 1$ is conserved over time $t \in [0, T]$.  \hfill \qed 

\smallskip

{\bf \noindent Proof of Lemma \ref{lem: sol of p. PME}.}
The proof is similar to that of Lemma \ref{lem: sol of PME}. Let  $ \bar{f} : [0, \infty) \to \R $ be a smooth function satisfying $ \bar{f}(u) = f(u)$ for $ \frac{1}{2\kappa} \le u \le \kappa + \frac{1}{2\kappa}$, $ \bar{f}'(u) > \epsilon^{-1} $ for $ 0 \le u \le \frac{1}{2\kappa} $ and $ \bar{f}'(u) < \epsilon $ for $ u \ge \kappa + \frac{1}{2\kappa}$, where $ \epsilon > 1 $ is some fixed constant.
Consider the PDE
\begin{align}	\label{eq: PME perturbed aux}
	\partial_t p^{\beta}(t,x) 
	= \diver \Big( \nabla \bar{f}\big(p^{\beta}(t,x)\big) + p^{\beta}(t,x) \nabla \beta(x)  \Big), \quad &\text{for } (t, x) \in (t_0, T] \times U
\end{align}
subject to the same initial and Neumann boundary conditions as in \eqref{eq: PME perturbed}. 
Again since $\bar{f}'(u) > \epsilon^{-1}$ for all $u \ge 0$, this PDE is uniformly parabolic. Therefore, standard quasilinear theory  implies that there exists a smooth solution $p^{\beta} \in C^{\infty}([t_0, T]$ $ \times \overline{U})$ to \eqref{eq: PME perturbed aux}. 

\smallskip

For a fixed  $\delta > 0$, we define 
\begin{align*}
	\lambda_{\delta} := \max_{t \in [t_0, \, t_0+\delta], \, x \in \overline{U}} \big\vert \partial_t p^{\beta}(t, x) \big\vert < \infty \quad \text{and} \quad \tau := \min\bigg(\delta, \, \frac{1}{2\kappa\lambda_{\delta}}   \bigg) > 0.
\end{align*}
 Let $T_{\beta} \coloneqq T \wedge (t_0 + \tau)$. Since $p^\beta_{t_0} = p_{t_0}$ by construction in \eqref{eq: PME perturbed},
\begin{equation*}
	\left| p^{\beta}(t, x) - p(t_0, x) \right|  = \left| p^{\beta}(t, x) - p^{\beta}(t_0, x) \right| \le  \int_{t_0}^{T_{\beta}} \big\vert \partial_t p^{\beta}(s, x) \big\vert \, \dd s \le \tau \lambda_{\delta} \le \frac{1}{2\kappa},
\end{equation*}
for every $(t,x) \in [t_0, T_\beta] \times \overline{U}$. From Lemma \ref{lem: sol of PME}, we have $ \kappa^{-1} \le p_{t_0} \le \kappa$, thus $ (2\kappa)^{-1} \le p^\beta_t\le \kappa + (2\kappa)^{-1}$ for all $t \in [t_0, T_\beta]$. This implies that $\bar{f}\big( p^{\beta}(t, x) \big) = f\big( p^{\beta}(t, x) \big)$ holds for every $(t,x) \in [t_0, T_\beta] \times \overline{U}$ and therefore $p^{\beta}$ also solves  \eqref{eq: PME perturbed}. 
	
Finally, for mass conservation, note that integration by parts gives us
\begin{align*}
	\frac{\dd}{\dd t} \int_{U} p^\beta(t,x) \, \dd x &= \int_{U} \diver \Big( \nabla f\big(p^{\beta}(t,x)\big) + p^{\beta}(t,x) \nabla \beta(x) \Big)  \,\dd x\\
	& = \int_{\partial U} \Big\langle \nabla f\big(p^{\beta}(t,x)\big) + p^{\beta}(t,x) \nabla \beta(x), \, n(x)\Big\rangle \, \dd x	= 0,
\end{align*}
where the last step follows from the no-flux boundary condition in \eqref{eq: PME perturbed} as well as Assumption \ref{assum.perturbation}.
\hfill \qed

\subsection{Proofs of Lemmas \ref{lem: stoch.rep} and \ref{lem: stoch.rep.p}} 

We will only prove Lemma \ref{lem: stoch.rep.p}, as the proof of Lemma \ref{lem: stoch.rep} is completely analogous.

For every $ (\tau, y) \in [t_0, T_\beta) \times \overline{U}$,  consider the SDE with reflection \eqref{eq: SDE.p} conditional on the initial position $ y $ at time $ \tau  $:
	\begin{equation} \label{eq: SDE X.tau.y}
		\begin{dcases}
			X^{\beta,\tau, y}_t = y - \int_{t_0}^t \nabla \beta(X^{\beta,\tau, y}_s) \, \dd s  + \int_{\tau}^t\sqrt{\frac{2f\big(p^\beta(s,X^{\beta,\tau, y}_s)\big)}{p^\beta(s,X^{\beta,\tau, y}_s)}}\dd W_s -  \int_{\tau}^t n(X^{\beta,\tau, y}_s)\, \dd L^{\beta, \tau,y}_s \in \overline{U},  \qquad  t \in [\tau, T_\beta],\\
			L^{\beta, \tau,y}_t = \int_\tau^t \mathrm{1}_{\{X^{\beta,\tau, y}_s \in \partial U\}} \, \dd L^{\beta, \tau,y}_s, \\
			L^{\beta, \tau,y}_\tau = 0 \text{ and } t \mapsto L^{\beta, \tau,y}_{t} \text{ is nondecreasing and continuous.}
		\end{dcases}
	\end{equation}
From Lemma \ref{lem: sol of p. PME} and Assumption \ref{assum1}(c), it is straightforward to check that diffusion coefficient is uniformly Lipschitz in the spatial variable. Thus by \cite[Theorem 3.1 and Remark 3.3]{lions1984stochastic}, the SDE with reflection \eqref{eq: SDE X.tau.y} has a pathwise unique, strong solution. Let $ \xi  $ be an independent $ \overline{U} $-valued random variable with distribution $ p_{t_0} $. Consider the process $ X^{\beta} $ given by $ X^{\beta}_{t_0} = \xi  $ and $ X^{\beta}_t = X^{\beta, t_0, \xi}_t $ for $ t \in (t_0, T_\beta] $. Similarly, let $ L^\beta $ be specified by $ L^\beta_{t_0} = 0$ and $ L^\beta_t = L^{\beta, t_0, \xi}_t $ for $ t \in (t_0, T_\beta] $.  Then $ (X^{\beta}, L^\beta)$ is the unique strong solution to \eqref{eq: SDE.p}. This completes the proof of the first part of the lemma.

Turning to the proof of the second part, we borrow  ideas from  \cite[Remark 3.1.2]{pilipenko2014introduction} and \cite[Chapter 5.7.B]{KS88}. Recall that $ p^\beta $ is fixed as the solution given in Lemma \ref{lem: sol of p. PME}. Consider the following \emph{backward Kolmogorov equation}:
\begin{align}	\label{eq: backward PDE}
	\begin{dcases}
		\partial_\tau q^\beta(\tau,y) + \frac{f(p^\beta(\tau,y))}{p^\beta(\tau,y)}\Delta_y q^\beta(\tau,y) - \left\langle \nabla_y  \beta(y), \nabla_y q^\beta(\tau ,y )  \right\rangle = 0, &\text{for } (\tau, y) \in (t_0,   T_\beta) \times U, \\
		\frac{\partial q^{\beta}(\tau, y)}{\partial n(y)} = 0, \qquad ~~~ &\text{for } y \in \partial U.
	\end{dcases}
\end{align}
It follows from \cite{ito1957boundary} that \eqref{eq: backward PDE} has a \emph{fundamental solution} $ G^\beta(\tau, y; t,x) $ defined for $ t_0 \le \tau < t \le T_\beta $ and $ x, y \in \overline{U} $. In particular, $G^\beta$ is nonnegative and for every $ \phi \in C(\overline{U})$ and $ t \in (\tau, T_\beta] $, the function
\begin{align} \label{eq: q(tau,y)}
	q^\beta(\tau,y) \coloneqq \int_{U} G^\beta(\tau,y;t,x) \, \phi(x) \,\dd x, 
\end{align}
 satisfies \eqref{eq: backward PDE} and the terminal condition
\begin{align}
	\lim_{\tau \uparrow t} q^\beta(\tau, y) = \phi(y), \quad \text{for all } y \in U.
\end{align}
If furthermore $\phi$ satisfies the no-flux boundary condition $\partial \phi \slash \partial n = 0$, then the above convergence holds uniformly in $\overline{U}$.

From the Feynman-Kac representation \cite[Theorem 3.1.1]{pilipenko2014introduction}, for any $ \phi \in C(\overline{U}) $,
	\begin{align} \label{eq: q(tau, y) 2}
		q^\beta(\tau,y) = \E\big[\phi(X^{\beta,\tau, y}_t)\big].
	\end{align}
Comparing \eqref{eq: q(tau,y)} with \eqref{eq: q(tau, y) 2}, we deduce that the transition probability density of $ X^{\beta,\tau, y}$ is given by  $ G^\beta $, i.e., 
\begin{align} \label{eq: transit prob den}
	\P\big(X^{\beta,\tau,y}_t \in A\big) = \int_A G^\beta(\tau, y; t,x) \, \dd x, \quad \text{for every Borel } A \subseteq \overline{U}.
\end{align}
For any fixed $y \in \overline{U}$, the function $ \Psi^\beta(t,x) \coloneqq G^\beta(t_0, y; t, x) $ satisfies the \emph{forward Kolmogorov equation}
\begin{align}
	\partial_t \Psi^\beta(t,x) = \text{div}_x \bigg(\nabla_x \Big(\frac{f\big(p^\beta(t,x)\big)}{p^\beta(t,x)} \Psi^\beta(t,x)\Big) + \Psi^\beta(t,x) \nabla_x  \beta(t,x)  \bigg),
\end{align}
which is the \emph{adjoint} of \eqref{eq: backward PDE}.
The probability density function of $ X^\beta_t $ with initial distribution $ p_{t_0} $ is then 
\begin{align}
	\widehat{p}^\beta(t,x) \coloneqq \int_{U} G^\beta(t_0,y;t,x) \, p_{t_0}(y) \, \dd y. 
\end{align} 
Together with \eqref{eq: assum noflux initial} of Assumption \ref{assum1}(c),  
we see that $p^\beta$ satisfies the linear uniformly parabolic PDE 
\begin{align*}
	\begin{dcases}
		\partial_t \widehat{p}^\beta(t,x) = \text{div} \bigg(\nabla \Big(\frac{f\big(p^\beta(t,x)\big)}{p^\beta(t,x)} \widehat{p}^\beta(t,x)\Big) + \widehat{p}^\beta(t,x) \nabla \beta(t,x)  \bigg),  &\text{for } (t,x) \in (t_0, T_\beta] \times U, \\ 
		\widehat{p}^\beta(t_0,x) = p(t_0,x), &\text{for }\ x \in \overline{U},\\
		\frac{\partial \widehat{p}^\beta(t,x)}{\partial n(x)} = 0, \qquad ~~~ &\text{for } x \in \partial U.
	\end{dcases}
\end{align*}
Since the solution to this PDE is unique, we deduce from Lemma \ref{lem: sol of p. PME} that $ \widehat{p}^\beta = p^\beta$. \hfill \qed

\subsection{Proof of Theorems \ref{thm: traj} and \ref{thm: traj.p.r.}} \label{subsec: proof.traj perturbed}
	We will only prove Theorem \ref{thm: traj.p.r.}, as similar arguments can be used to show Theorem \ref{thm: traj}. We first prove \eqref{eq: semi-mart decom.p}. Recall the function $ \varphi $ defined in \eqref{eq: D process} and note for later use the simple identities
	\begin{align} \label{eq: simple.iden}
		\varphi(u) = h(u) -\frac{f(u)}{u}, \quad h(u) = \varphi'(u)u + \varphi(u), \quad \Phi''(u) = \varphi''(u)u + 2 \varphi'(u) \quad \text{for all } u > 0.
	\end{align}
 By writing $v^{\beta}(t, x) = \varphi\big(p^\beta(t, x)\big)$, we deduce from \eqref{eq: PME perturbed} that
	\begin{align} \label{eq: v PPDE}
		\partial_t v^{\beta}(t, x) = \varphi'\big(p^{\beta}(t, x)\big) \partial_t p^{\beta} (t, x) 
		= \varphi'\big(p^{\beta}(t, x) \big) \, \diver \Big( \nabla f\big(p^{\beta}(t, x)\big) + p^{\beta}(t, x) \nabla \beta(x)  \Big). 
	\end{align}
	Using It\^{o}'s lemma along with \eqref{eq: SDE.p} and \eqref{eq: v PPDE}, we see that the dynamics of the perturbed entropy process satisfies
	\begin{align*}
		\dd v^{\beta} (t, X^{\beta}_t) &= \bigg(\varphi'(p^{\beta}) \, \diver \big( \nabla f(p^{\beta}) + p^{\beta} \nabla \beta \big) +  \frac{f(p^{\beta})}{p^\beta} \Delta v^{\beta} - \langle \nabla v^{\beta}, \, \nabla \beta \rangle\bigg) (t, X^{\beta}_t) \, \dd t  \\
		&\qquad  \qquad  + \bigg\langle \bigg(\sqrt{\frac{2f(p^{\beta})}{p^\beta}} \nabla v^{\beta}\bigg) (t, X^{\beta}_t), \, \dd W_t  \bigg\rangle
		- \big\langle \nabla v^{\beta}  (t, X^{\beta}_t), \, n(X^{\beta}_t)  \big\rangle \, \dd L^\beta_t 
		\\
		&= D^\beta(t, X^{\beta}_t) \, \dd t + M^{\beta}_t - \varphi'\big(p^{\beta}(t, X_t)\big) \big\langle \nabla p^{\beta}(t, X^{\beta}_t), \, n(X^{\beta}_t) \big\rangle \, \dd L^\beta_t.
	\end{align*}
	Note that the last line in \eqref{eq: SDE.intro} ensures that the reflecting term $L^\beta$ only increases when $X^{\beta}$ is on the boundary. In conjunction with the no-flux boundary condition in \eqref{eq: PME perturbed}, we see that the last term above is zero. This completes the proof of \eqref{eq: semi-mart decom.p}.
	
	Next, to see that the local martingale $ M^\beta$ in \eqref{eq: cum en diss process.p} is in fact a true $ L^2 $-martingale, note that the quadratic variation of $ M^\beta $ is given by 
	\begin{align*}
	\E \Big[\langle M^{\beta}, \, M^{\beta}\rangle_{T_{\beta}} \Big] = \E \int_{t_0}^{T_{\beta}} \bigg( \frac{2f(p^{\beta})}{p^\beta} \vert\nabla v^{\beta}\vert^2 \bigg) (t, X^{\beta}_t) \, \dd t = \E \int_{t_0}^{T_{\beta}} \bigg( \frac{2f(p^{\beta})^3 }{(p^\beta)^5}\vert\nabla p^{\beta}\vert^2 \bigg) (t, X^{\beta}_t) \, \dd t.
    \end{align*}
	We claim that the above quantity is finite. Indeed, due to the properties of the solution $p^\beta$ in Lemma \ref{lem: sol of p. PME}, 
	the above expectation is bounded by 
	\begin{align*}
		2(T_\beta - t_0) (2\kappa)^5 f\left(\frac{3}{2\kappa}\right)^3 	\max_{[t_0, T_{\beta}] \times \overline{U}} \vert \nabla p^\beta(t, x) \vert^2 < \infty.
	\end{align*}
    Therefore, it follows from \cite[Corollary IV.1.25]{RY99} that $ \big(M^\beta_t\big)$ is an $ L^2 $-martingale.

    Finally, in order to show \eqref{eq: fish.inf.exp.p}, we take expectation in the first equation in \eqref{eq: cum en diss process.p} and use Fubini's theorem to get 
	\begin{align} \label{eq: sum expect}
		\E [F^{\beta}_t ] 
		= \int_{t_0}^t \E \big[D^{\beta}(s, X^{\beta}_s)\big] \, \dd s = \sum_{i=1}^3\int_{t_0}^t \E \big[D^{\beta}_i(s, X^{\beta}_s)\big] \, \dd s
	\end{align}
	where 
	\begin{align*}
		D^{\beta}_1 \coloneqq \varphi'(p^{\beta}) \, \diver \Big(\nabla f\big(p^{\beta}\big) + p^{\beta} \nabla \beta  \Big), \quad
		D^{\beta}_2  \coloneqq  \frac{f(p^{\beta})}{p^\beta} \Delta v^{\beta}, \quad \text{and} 
		\quad D^{\beta}_3  \coloneqq - \big\langle \nabla v^{\beta}, \, \nabla \beta \big\rangle.
	\end{align*}
    We now evaluate each of the expectations in \eqref{eq: sum expect}. Integrating by parts, we have
	\begin{align} \label{eq: E[D_1]}
	\begin{split}
	    \E \big[D^{\beta}_1(t, X^{\beta}_t)\big] &= \int_{U} \varphi'(p^{\beta}) p^{\beta} \, \diver \Big( \nabla f\big(p^{\beta}\big) + p^{\beta} \nabla \beta  \Big) (t, x) \,  \dd x \\
		&= - \int_{U} \Big\langle \nabla \Big(\varphi'(p^{\beta})p^{\beta}\Big), \nabla f\big(p^{\beta}\big) + p^{\beta} \nabla \beta \Big\rangle (t, x)\, \dd x + C,
	\end{split}
	\end{align}
	where $C$ is the boundary term given by
	\begin{align*}
	    C \coloneqq \int_{\partial U} \varphi'(p^\beta) f'(p^\beta) p^\beta \Big\langle \nabla p^\beta + p^\beta \nabla \beta ,n \Big\rangle (t,x) \, \dd x.
	\end{align*}
	From the no-flux boundary condition in \eqref{eq: PME perturbed} and Assumption \ref{assum.perturbation}, we see that $C = 0$. Similarly, we see that $\E \big[D_2^{\beta}(t, X^{\beta}_t)\big] $ and $\E \big[D_3^{\beta}(t, X^{\beta}_t)\big] $ are respectively equal to  
	\begin{align} \label{eq: E[D_2]}
		 - \int_{U} \big\langle \nabla f(p^{\beta}), \nabla v^{\beta} \big\rangle (t, x) \, \dd x, \quad  \text{and} \quad - \int_{U} \big\langle \nabla v^\beta, \nabla \beta \big\rangle (t,x) \, \dd x.
	\end{align}
	Assembling them gives  
	\begin{align} \label{eq: D expectation}
	\begin{split}
			\E [D^{\beta}(t, X^{\beta}_t)] &= - \int_{U} \Big\langle \nabla \big(\varphi'(p^{\beta})p^{\beta} + v^{\beta} \big), \nabla f(p^{\beta})\Big\rangle (t, x) \, \dd x \\
		&\qquad \qquad- \int_{U}\Big\langle \nabla \big(\varphi'(p^{\beta})p^{\beta}\big) + \nabla v^{\beta}, \, p^{\beta} \nabla \beta  \Big\rangle (t, x) \, \dd x.
	\end{split}
	\end{align}
	Using the third identity in \eqref{eq: simple.iden}, we see that the first integrand above is equal to 
	\begin{align*}
		\big(\varphi''(p^{\beta})p^{\beta} + 2 \varphi'(p^{\beta})\big) f'(p^{\beta}) \vert \nabla p^{\beta}\vert^2 = \Phi''(p^{\beta}) f'(p^{\beta}) \vert \nabla p^{\beta}\vert^2,
	\end{align*} so the first integral in \eqref{eq: D expectation} is
	\begin{align*}
		 -\int_{U} \Big( \big\vert \Phi''(p^{\beta})\nabla p^{\beta} \big\vert^2 p^{\beta} \Big) (t, x) \, \dd x 
		 = -I\big( p^{\beta}_t\big) > - \infty,
	\end{align*}
	where the last step follows from the boundedness of $p^{\beta}$ and $\nabla p^{\beta}$ implied by Lemma \ref{lem: sol of p. PME}. Similarly, using the second identity in \eqref{eq: simple.iden}, we see that the second integral in \eqref{eq: D expectation} is equal to
	\begin{equation*}
		-\E \Big[\big\langle \nabla \big(h(p^{\beta})\big), \, \nabla \beta \big\rangle(t, X^{\beta}_t)\Big] > -\infty.
	\end{equation*}
Putting them together completes the proof of \eqref{eq: fish.inf.exp.p}. \hfill \qed

\subsection{Proofs of Corollaries \ref{cor: deBruijn} and \ref{cor: deBruijn.p}}
	We will only prove Corollary \ref{cor: deBruijn.p}, as the proof of Corollary \ref{cor: deBruijn}  proceeds in the same way.
	
	Taking expectation in \eqref{eq: semi-mart decom.p} and using the martingale property of $ M^\beta $ in  \eqref{eq: cum en diss process.p} as well as \eqref{eq: fish.inf.exp.p}, we have
	\begin{align}
		\FFF\big(p^{\beta}_t\big) - \FFF\big(p^{\beta}_{t_0}\big)
		&= \E\Big[ v^{\beta}(t, X^{\beta}_{t})-v^{\beta}(t_0, X^{\beta}_{t_0})\Big] \nonumber\\
		&
		=-\int_{t_0}^{t} I( p^{\beta}_s)  \, \dd s  - \int_{t_0}^t  \E \Big[ \big\langle \nabla h \big(p^{\beta}(s, X^{\beta}_s) \big), \, \nabla \beta(X^{\beta}_s) \big\rangle \Big] \, \dd s,  \label{eq: second integrand}
	\end{align}
	which proves \eqref{eq: entro diss integral.p}.
	
	Turning to the proof of \eqref{eq: entro diss iden.p}, note that from \eqref{eq: h, Phi def} we have
	\begin{align} \label{eq: cont s}
		\E \Big[ \big\langle \nabla h \big(p^{\beta}(s, X^{\beta}_s) \big), \, \nabla \beta(X^{\beta}_s) \big\rangle \Big] = \int_{U} f'\big(p^\beta(s,x)\big) \big\langle \nabla p^\beta(s,x), \nabla \beta(x) \big\rangle \, \dd x.
	\end{align}
	From the continuity of $(s ,x) \mapsto f'\big(p^\beta(s,x)\big) \big\langle \nabla p^\beta(s,x), \nabla \beta(x) \big\rangle $, we see that the expression of \eqref{eq: cont s} is continuous as a function of $ s $, thus
	\begin{equation}
		\frac{\mathrm{d}}{\mathrm{d}t} \bigg\vert_{t=t_0^+} \int_{t_0}^t  \E\Big[ \big\langle \nabla h \big(p^{\beta}(s, X^{\beta}_s) \big), \, \nabla \beta(X^{\beta}_s) \big\rangle \Big] \dd s = \E \Big[ \big\langle \nabla h \big(p^{\beta}(t_0, X_{t_0}) \big), \, \nabla \beta(X_{t_0}) \big\rangle \Big],
	\end{equation}
	where the last equality is due to the fact that $ X^\beta_{t_0} = X_{t_0}$ by construction in \eqref{eq: PME perturbed}.  Similarly, from the continuity of $ t \mapsto I\big(p^\beta_t\big) $, we have
	\begin{equation}
		\frac{\mathrm{d}}{\mathrm{d}t} \bigg\vert_{t=t_0^+} \int_{t_0}^t I\big(p^{\beta}_u\big) \, \mathrm{d}u
		= I\big(p^{\beta}_{t_0}\big)
		= I(p_{t_0}),
	\end{equation}
	and the identity \eqref{eq: entro diss iden.p} follows.
	
	Finally, to show \eqref{eq: conditional entro diss perturb 1},  the martingale property of $M^\beta$ implies that the numerator on the left-hand side of \eqref{eq: conditional entro diss perturb 1} is equal to
	\begin{equation*}
		\mathbb{E} \Big[ F^\beta_{t} - F^\beta_{t_0} \, \big| \, \FF_{t_0} \Big]
		= \mathbb{E} \left[ \int_{t_0}^{t} D^\beta (u, X^{\beta}_u) \,\dd u \, \Big| \, \FF_{t_0} \right].
	\end{equation*}
	From the continuity of $u \mapsto D^\beta(u, X_u)$, we have
	\begin{equation*}
		\lim_{t \downarrow t_0} \frac{1}{t-t_0} \int_{t_0}^{t} D^\beta(u, X^{\beta}_u) \, \dd u = D^\beta(t_0, X_{t_0}), \quad \text{a.s.}.
	\end{equation*}
	Moreover, the properties of $p^\beta$ from Lemma \ref{lem: sol of p. PME} implies that $D^\beta$ is uniformly bounded on $[t_0, T_\beta] \times \overline{U}$. Therefore, by the bounded convergence theorem,
	\begin{equation*}
		\lim_{t \downarrow t_0} \mathbb{E} \left[ \frac{1}{t-t_0} \int_{t_0}^{t} D^\beta(u, X^{\beta}_u) \dd u \right] 
		= \mathbb{E} \big[D^\beta(t_0, X_{t_0})\big].
	\end{equation*}
	We now apply Scheff\'{e}'s lemma to obtain
	\begin{equation*}
		\lim_{t \downarrow t_0} \left\Vert \frac{1}{t-t_0} \int_{t_0}^{t} D^\beta(u, X^{\beta}_u) \dd u - D^\beta(t_0, X_{t_0}) \right\Vert_{L^1} = 0.
	\end{equation*}
	To complete the proof, use Jensen's inequality and the tower property to get
	\begin{align*}
		\bigg\Vert \mathbb{E} \bigg[ \frac{1}{t-t_0} \int_{t_0}^{t} D^\beta(u, X^{\beta}_u) \, \dd u \, \Big| \, \FF_{t_0} \bigg] &- D^\beta(t_0, X_{t_0}) \bigg\Vert_{L^1} \\
		&\le \left\Vert \frac{1}{t-t_0} \int_{t_0}^{t} D^\beta(u, X^{\beta}_u) \, \dd u - D^\beta(t_0, X_{t_0}) \right\Vert_{L^1},
	\end{align*}
	and the $L^1$-convergence in \eqref{eq: conditional entro diss perturb 1} follows. 
	 \hfill \qed
	
\medskip

\subsection{Proof of Lemma~\ref{thm: wass}} \label{sec: Wasserstein}
Since the proofs of \eqref{eq: wass.limit} and \eqref{eq: wass.limit.perturb} are very similar, we will only prove \eqref{eq: wass.limit.perturb}.
We first rewrite the PDE in \eqref{eq: PME perturbed} as a continuity equation
\begin{align}	\label{eq: cont.eq}
	\partial_t p^{\beta}(t,x) 
	+ \diver \Big(p^{\beta}(t,x) \, \vel^\beta(t,x)\Big) = 0, \quad \text{for } (t, x) \in (t_0, T_\beta) \times U,
\end{align}
where $\vel^\beta : [t_0, T_\beta] \times \overline{U} \to \R^d $ is the \emph{velocity field}, defined by
\begin{align} \label{eq: v.fields.beta}
	\vel^\beta(t,x) \coloneqq -  \nabla \left[ \beta+ h\big(p^\beta\big) \right] (t,x).
\end{align}
We see from Assumption \ref{assum.perturbation} and Lemma \ref{lem: sol of p. PME} that $\vel^{\beta}(t_0, \cdot)$ is the gradient of a smooth function. 
For each $ x \in \overline{U}$, consider the curved flow $ \Lambda^\beta_t $ associated with $\vel^{\beta}$, specified by
\begin{align} \label{eq: curved flow}
	\frac{\dd}{\dd t} \Lambda^{\beta}_t (x)= \vel^{\beta}\big(t,\Lambda^{\beta}_t(x) \big), \quad t \in [t_0, T_\beta], \quad \Lambda^{\beta}_{t_0}(x) = x.
\end{align}
By the Cauchy-Lipschitz theorem, there exists a unique solution $ t \mapsto \Lambda^\beta_t(x) \in \overline{U} $ to \eqref{eq: curved flow}. Moreover, it follows from \cite[Theorem 5.34]{Vil03} that $\Lambda^{\beta}_t$ pushes forward $p^{\beta}_{t_0}$ to $p^{\beta}_t$, in the sense that $p^\beta_{t_0} \circ (\Lambda^\beta_t)^{-1} = p^{\beta}_t$. Note also that the Jacobian of $\Lambda^\beta_t$ is given by
\begin{align*}
	\nabla \Lambda^\beta_t(x) = I - \int_{t_0}^t \nabla^2\left[\big(\beta + h(p^\beta)\big) \big(s, \Lambda^\beta_s(x)\big)\right] 	\dd s.
\end{align*}
Therefore, by setting 
\begin{align*}
    K := \max\left\{\Big|\partial_{ij} \big(\beta + h(p^\beta)\big)(t, x) \Big|: t \in [t_0, T_\beta], x \in \overline{U}, i, j =1, \dots, n\right\},
\end{align*}
 we see that for any $ t \in [t_0, t_0 + K^{-1})$, the Jacobian $ \nabla \Lambda_t^\beta (x)$ is positive-semidefinite for all $ x \in \overline{U} $, so $ \Lambda^\beta_t $ is the gradient of a convex function. Hence, by Brenier's theorem \cite[Theorem 2.12(ii)]{Vil03}, $ \Lambda^\beta_t $ is the optimal transport map from $ p^\beta_{t_0} $ to $ p^\beta_t $, i.e., 
\begin{align*}
	\W^2_2\big(p^\beta_{t_0}, p^\beta_t\big) = \E \left[\big| \Lambda^\beta_t(X^{\beta}_{t_0}) - X_{t_0} \big|^2\right] = \E \left[\bigg| \int_{t_0}^t \vel^\beta\big(s, \Lambda^\beta_s(X^{\beta}_{t_0})\big) \, \dd s\bigg|^2\right].  
\end{align*}
Now, the continuity of $ t \mapsto \vel^\beta\big(t, \Lambda^\beta_t(x)\big) $ implies
\begin{align} \label{eq: a.s.conv}
	\left|\frac{1}{t - t_0} \int_{t_0}^t \vel^\beta \big(s, \Lambda^\beta_s(X^{\beta}_{t_0})\big) \, \dd s \right|^2\stackrel{t \downarrow t_0}{\longrightarrow} \big|\vel^\beta(t_0, X_{t_0})\big|^2, \quad \text{a.s}.
\end{align}
Moreover, by Jensen's inequality, the random variable on the left-hand side above is bounded by
\begin{align*}
	\frac{1}{t -t_0} \int_{t_0}^t \big|\vel^\beta(s, \Lambda^\beta_s(X^{\beta}_{t_0})) \big|^2 \dd s \le \max\big\{ \vert \vel^\beta(t,x) \vert^2 : t \in [t_0, T_\beta], x \in \overline{U}\big\} < \infty,
\end{align*}
where the last step follows from the aforementioned fact that $ \vel^\beta (t, \cdot)$  is the gradient of a smooth function. Consequently, by the bounded convergence theorem,
\begin{align}
	\lim_{t \downarrow t_0} \frac{\W_2(p_{t_0}^\beta, p^\beta_t)}{t - t_0} = \lim_{t \downarrow t_0} \bigg\Vert  \frac{1}{t - t_0}\int_{t_0}^t \vel^\beta\big(s, \Lambda^\beta_s(X^{\beta}_{t_0})\big) \, \dd s\bigg\Vert _{L^2} = \big\Vert\vel^\beta(t_0, X_{t_0})\big\Vert_{L^2},
\end{align}
where the last step follows from the fact that $X^\beta_{t_0} = X_{t_0}$ by construction in \eqref{eq: PME perturbed}. Recalling the expression of $ \vel^\beta $ in \eqref{eq: v.fields.beta}, we arrive at \eqref{eq: wass.limit.perturb}.
 \hfill\qedsymbol

\subsection{Proof of Theorem~\ref{thm: steepest.decsent}}
The identity \eqref{eq: FW} follows from \eqref{eq: entro diss iden2} of Corollary \ref{cor: deBruijn} and \eqref{eq: wass.limit} of Lemma \ref{thm: wass}. Similarly, for \eqref{eq: FW.p}, we deduce from \eqref{eq: entro diss iden.p} of Corollary \ref{cor: deBruijn.p} and \eqref{eq: wass.limit.perturb} of Lemma \ref{thm: wass} that
\begin{align} \label{eq: wass.lim.pf}
	  |\partial \mathscr{F}|_{\W_2} (p^\beta_{t_0})
	&=  - \frac{I(p_{t_0}) + \E \Big[ \big\langle \nabla h \big(p(t_0, X_{t_0}) \big), \, \nabla \beta(X_{t_0})\Big] }{\Big\Vert\nabla h\big(p(t_0, X_{t_0})\big) + \nabla \beta(X_{t_0})  \Big\Vert_{L^2}}.
\end{align}
Moreover, we see from \eqref{eq: energy diss func} and \eqref{eq: h, Phi def} that the entropy dissipation functional can be expressed as
\begin{align*}
	I(p_{t_0}) =  \E\Big[\big\vert \nabla h \big(p(t_0, X_{t_0})\big)\big\vert^2\Big].
\end{align*}
Putting this back into \eqref{eq: wass.lim.pf} yields \eqref{eq: FW.p}. Finally, the inequality \eqref{eq: FW-FW.p}  follows from the Cauchy-Schwarz inequality. \hfill \qed

\subsection{Proof of Proposition \ref{prop: FWI}}
For $t \in [0,1)$, let us denote by
 \begin{align} \label{eq: def T_t}
     T_t \coloneqq (1 - t) \textrm{Id} + t \nabla \psi
\end{align}
the optimal transport map from $\rho_0$ to $\rho_t$. 
Note that $T_t$ is injective by \cite[Theorem 5.49]{Vil03}, so its inverse exists. By \cite[Theorem 5.51(ii)]{Vil03}, the probability density functions $(\rho_t)_{t \in (0,1)}$  satisfy the continuity equation
\begin{equation} \label{eq: contunity equation, HWI}
    \partial_t \rho_t(x) + \textrm{div} \big(\rho_t(x) \, \velo_t(x)\big) = 0,
\end{equation}
where $\velo: [0, 1) \times \overline{U} \to \R^d $ is the velocity field defined by
\begin{equation} \label{eq: v.f. FWI}
    \velo_t(x) \coloneqq (\nabla \psi - \textrm{Id}) \circ (T_t)^{-1}(x).
\end{equation}
In conjunction with \eqref{eq: def T_t}, we see that $T_t$ satisfies the integral equation
\begin{equation} \label{eq: T_t flow}
    T_t(x) = x + \int_0^t \velo_s \big(T_s(x)\big) \dd s.
\end{equation}

We now switch to probabilistic notations. On a sufficiently rich probability space, let $Z_0$ be a random variable with distribution $\rho_0$ and let $Z_t \coloneqq T_t(Z_0)$ for $t \in (0, 1)$. On account of \eqref{eq: T_t flow}, we have 
\begin{equation} \label{eq: Z dym}
    Z_t = Z_0 + \int_0^t \velo_s(Z_s)\,\dd s.
\end{equation}
Together with \eqref{eq: contunity equation, HWI}, we deduce
\begin{equation*}
    \dd \rho_t(Z_t) = \partial_t\rho_t(Z_t) \, \dd t + \big\langle \nabla \rho_t(Z_t), \dd Z_t\big\rangle = -\rho_t(Z_t) \, \textrm{div}\big(\velo_t(Z_t)\big) \dd t.
\end{equation*}
Recalling the function $\varphi $ in \eqref{eq: D process}, we have
\begin{equation*}
    \dd \varphi\big(\rho_t(Z_t)\big) = -\varphi'\big(\rho_t(Z_t)\big) \rho_t(Z_t) \, \textrm{div}\big(\velo_t(Z_t)\big) \dd t = -\frac{f\big(\rho_t(Z_t)\big)}{\rho_t(Z_t)} \textrm{div}\big(\velo_t(Z_t)\big) \dd t,
\end{equation*}
where the last step follows from the identity $f(u) = u h(u) - \Phi(u)$, valid for all $u \ge 0$.
Integrating from $0$ to $t$ and taking expectation yield
\begin{align}
    \FFF(\rho_t) - \FFF(\rho_0) & = - \int_0^t \int_{U} f\big(\rho_s(z)\big) \textrm{div}\big(\velo_s(z)\big)  \dd z \dd s \nonumber\\
    &= \int_0^t \int_{U}  \big\langle \nabla f\big(\rho_s(z)\big), \velo_s(z) \big\rangle \, \dd z \dd s - \int_0^t \int_{\partial U} f \big(\rho_s(z)\big) \big\langle \velo_s(z), n(z) \big\rangle \, \dd z \, \dd s. \label{eq: FWI ent.diss}
\end{align}
It follows from Assumption \ref{FWI assum}(b) that $ \nabla \psi (x) = x $ for all $ x \in \partial U $.  Therefore, we see from \eqref{eq: v.f. FWI} that $ \velo_t(x) = 0 $ for all $ x \in \partial U $. Hence, the last integral in \eqref{eq: FWI ent.diss} vanishes. Letting $t \downarrow 0$ yields
\begin{align*}
    \frac{\dd}{\dd t} \Big|^{+}_{t =0} \FFF(\rho_t) &=  \int_{U}  \big\langle \nabla f \big(\rho_0(z)\big), \velo_0(z) \big\rangle \, \dd z 
\end{align*}
Recalling the definition of $ \velo_0 $ in \eqref{eq: v.f. FWI} completes the proof. \hfill \qed

\subsection{Proof of Theorem \ref{thm: FWI}}
   From \cite[Theorem 5.15(i)]{Vil03}, Assumption \ref{FWI assum}(c) implies that the entropy functional $ \mathscr{F} $ is \emph{displacement convex}. In other words, 
\begin{equation*} 
    \frac{\dd^2}{\dd t^2} \FFF(\rho_t) \ge 0, \quad \text{for} \quad t \in [0,1].
\end{equation*}
Therefore,  Taylor's theorem and Proposition \ref{prop: FWI} give us
\begin{align*}
    \FFF(\rho_1) &= \FFF(\rho_0) + \frac{\dd}{\dd t} \Big|^{+}_{t =0} \FFF(\rho_t) + \int_0^1 (1-t) \frac{\dd^2}{\dd t^2} \FFF(\rho_t) \, \dd t \\
    & \ge \FFF(\rho_0) + \int_{U} \big\langle \nabla f\big(\rho_0(z)\big), \nabla \psi(z) - z \big\rangle \, \dd z,
\end{align*}
which proves the first inequality in \eqref{eq: first.func.ineq}. The second inequality is a simple consequence of the  
 Cauchy-Schwarz inequality; the expression in the middle of  \eqref{eq: first.func.ineq} is bounded from above by
	\begin{align*}
		 \sqrt{\int_{U} \frac{\big|\nabla f \big(\rho_0(z) \big)\big|^2 }{\rho_0(z)}\dd z} \, \sqrt{\int_{U} \big|\nabla \psi(z) - z \big|^2 \rho_0(z) dz} = \sqrt{I(\rho_0)} \, \W_2(\rho_0, \rho_1).
	\end{align*}
\hfill \qed

\section*{Acknowledgements}
We are indebted to Ioannis Karatzas for suggesting this problem and for many helpful discussions. L.C. Yeung is partially supported by the National Science Foundation (NSF) under grant NSF-DMS-20-04997.

\bibliography{A_trajectorial_approach_to_entropy_dissipation_for_degenerate_parabolic_equations}
\bibliographystyle{amsplain}

\end{document}